\documentclass[a4paper,11pt]{amsart}
\usepackage[utf8]{inputenc}
\usepackage{CJKutf8}
\usepackage{amsmath,amsthm,amssymb}
\usepackage{graphicx}
\usepackage{subcaption}
\usepackage{float}
\usepackage{cancel}
\usepackage{slashed}
\usepackage{comment}
\newcommand{\mychoice}[3]{#1
}
\mychoice{
\newcommand{\plabel}[1]{ \label{#1}}
\newcommand{\gbibitem}[1]{ \bibitem{#1}}
\newcommand{\snewpage}{}

\specialcomment{commentx}{}{}\excludecomment{commentx}
\specialcomment{commenty}{}{}\excludecomment{commenty}
\newcommand{\rechoicecomm}[1]{}
}{
\usepackage{xcolor}
\newcommand{\plabel}[1]{ \label{#1}\rlap{\smash{${}^{^{[#1]}}$}}}
\newcommand{\gbibitem}[1]{ \bibitem{#1}\rlap{\smash{${}^{^{[#1]}}$}}}
\newcommand{\snewpage}{\newpage}

\newenvironment{commentx}{\color{magenta} }{\color{black} }
\newenvironment{commenty}{\color{blue} }{\color{black} }
\newcommand{\rechoicecomm}[1]{#1}
}{
\usepackage{xcolor}
\newcommand{\plabel}[1]{ \label{#1}}
\newcommand{\gbibitem}[1]{ \bibitem{#1}}
\newcommand{\snewpage}{}

\specialcomment{commentx}{}{}\excludecomment{commentx}
\newenvironment{commenty}{\scriptsize }{\normalsize }
\newcommand{\rechoicecomm}[1]{}
}
\DeclareMathOperator{\sgn}{sgn}
\DeclareMathOperator{\Id}{Id}

\DeclareMathOperator{\des}{des}

\DeclareMathOperator{\artanh}{artanh}

\DeclareMathOperator{\asc}{asc}

\DeclareMathOperator{\Dbar}{D}

\DeclareMathOperator{\Rexp}{exp_{R}}

\theoremstyle{definition}
\newtheorem{point}{}[section]

\newtheorem{remark}[point]{Remark}
\newtheorem{example}[point]{Example}

\theoremstyle{plain}

\newtheorem{lemma}[point]{Lemma}

\newtheorem{theorem}[point]{Theorem}
\theoremstyle{remark}

\newcommand{\leaveout}[1]{}

\newcommand{\qedexer}{\renewcommand{\qedsymbol}{$\diamondsuit$}\qed\renewcommand{\qedsymbol}{$\Box$}}
\newcommand{\qedremark}{\renewcommand{\qedsymbol}{$\triangle$}\qed\renewcommand{\qedsymbol}{$\Box$}}

\newcommand{\eqedexer}{
\renewcommand{\qedsymbol}{$\diamondsuit$}
\pushQED{\qed}
\qedhere
\popQED
\renewcommand{\qedsymbol}{$\Box$}
}
\newcommand{\eqedremark}{
\renewcommand{\qedsymbol}{$\triangle$}
\pushQED{\qed}
\qedhere
\popQED
\renewcommand{\qedsymbol}{$\Box$}
}

\newcommand{\marginextend}[1]{ \addtolength{\oddsidemargin}{-#1}  \addtolength{\evensidemargin}{-#1}
  \addtolength{\textwidth}{#1}\addtolength{\textwidth}{#1}}
\newcommand{\updownextend}[1]{ \addtolength{\topmargin}{-#1}  \addtolength{\textheight}{#1}
\addtolength{\textheight}{#1}}
\marginextend{1cm}
\updownextend{0cm}
\allowdisplaybreaks[4]
\title[Convergence estimates for the Magnus expansion IE]{Convergence estimates for the Magnus expansion IE.
Finite dimensional Banach algebras}
\author{Gyula Lakos}
\address{Alfréd Rényi Institute of Mathematics, Reáltanoda utca 13-15, Budapest, H--1053, Hungary}
\email{lakos@renyi.hu}
\keywords{Magnus expansion, Baker--Campbell--Hausdorff expansion, convergence estimates, resolvent method}
\subjclass[2020]{Primary: 46H30, Secondary:  16W80.}
\begin{document}
\begin{abstract}
We  review and provide simplified proofs related to the Magnus expansion, and improve convergence estimates.
Observations and improvements concerning the Baker--Campbell--Hausdorff expansion are also made.

In this Part IE, we consider the case of finite dimensional Banach algebras.
We show that Magnus expansion is convergent (and works in logarithmic sense) if the cumulative norm
$< 2+\varepsilon_{n}$, where $\varepsilon_{n}$ is a positive number depending on the
dimension $n$ of the Banach algebra.
We also show concrete finite-dimensional counterexamples of multiple Baker--Campbell--Hausdorff type
for any cumulative norm $>2$ (necessarily of possibly great dimension).
\end{abstract}
\maketitle
\section*{Introduction}
This paper is a direct continuation of Part I \cite{L1}.
Here we will consider finite dimensional Banach algebras but which  are
 not far from the extremal case of Part I.

\textbf{Introduction to the setting of finite dimensional Banach algebras.}
Since its inception, Magnus \cite{M},
 finite dimensional normed algebras, which are automatically Banach algebras,
 play a prominent role in the study and applications of the Magnus expansion,
 cf.~Blanes, Casas, Oteo, Ros \cite{BCOR}.
In the case of operators on Hilbert spaces, finitely dimensionality is not a strong restriction
 with respect to the convergence radius in terms of the cumulative norm:
The counterexamples by Sch\"affer \cite{Scha}  (indirectly), Vinokurov \cite{V}, Moan \cite{Ma},
 and the convergence results of Moan, Niesen \cite{MN}, Casas \cite{Ca} show
 that the critical case for the convergence of the Magnus expansion for Hilbert space operators
 is the cumulative norm $\pi$, and it is represented by concrete counterexamples over $2\times2$ real matrices.
(See Part II  \cite{L2} for a detailed discussion.)
In the case of general, possibly infinite dimensional  Banach algebras, results of Moan, Oteo \cite{MO} and Part I \cite{L1}
 show that the critical convergence radius in terms of the cumulative norm is $2$,
 and counterexamples for convergence with cumulative norm $2$ must need to be somewhat wild.
It is not immediately clear that how finite dimensionality affects these results; but, as it turns out, positively.
An ease in the finite dimensional case that issues of convergence do not depend on the choice of the norm, as all norms are equivalent.
Selecting among the convergent or divergent measures in terms of the possible norms is a different matter, however.

\snewpage
\textbf{Outline of content.}
In Section \ref{sec:DivFin}, we start with considering the ``Minimal Examples'',
 which are counterexamples to the convergence of the BCH expansion in the setting of $2\times2$
 upper triangular matrix algebras; these are of cumulative norm $\pi+\varepsilon$.
Next, we exhibit examples of multiple Baker--Campbell--Hausdorff type for divergence
 with cumulative norm $2\frac{n}{n-1}$ in some $n\times n$ matrix algebras (which are, of course, of dimension  $d=n^2$), where $n\geq2$.
Although examples are sometimes results of lucky choices, here will explain how these latter finite dimensional
 counterexamples relate to the infinite dimensional ones of Part I \cite{L1}.
In Section \ref{sec:ConFin},
 we show that Magnus expansion is convergent  and works in logarithmic sense  if the cumulative norm
 $< 2+\varepsilon_{n}$, where $\varepsilon_{n}$ is a positive number depending on the (real)
 dimension $n$ of the Banach algebra.
This section draws heavily on the resolvent method introduced in the formal case by
 Mielnik,  Pleba\'{n}ski \cite{MP}, and developed further in Part I \cite{L1}.
Interestingly, methods of discrete geometry are used; notably,
 we use  the estimates of Rogers \cite{R} on packing densities of centrally convex bodies.

\textbf{Acknowledgements.}
The author is grateful to Fernando Casas and Kurusch Eb\-ra\-hi\-mi-Fard for organizing venues for furthering research,
 and also to Pierre-Louis Giscard, Ilya Kuprov, and Stefano Pozza for inspiring conversations.

%\snewpage

\begin{commentx}
\tableofcontents
\end{commentx}

\snewpage

\section{Divergence for finite dimensional Banach algebras}
\plabel{sec:DivFin}
Here we will see how the higher dimensionality leads to divergent Magnus series with   measures of cumulative norm
 approaching $2$ from above.

\subsection{Low dimensional counterexamples}
\plabel{ss:low}
~\\

Real and complex algebras of low dimension can be classified, cf. Study \cite{Stu}, Scheffers \cite{Sche},
 Kobayashi, Shirayanagi,  Takahasi,  Tsukada \cite{KSTT}.
It is trivial that if the dimension $d$ of the unital algebra $\mathfrak A$ is $1$ or $2$, then the algebra is commutative.
Recall that in the commutative case the Magnus series is trivially convergent;
 the logarithmic version of Magnus formula holds if the cumulative norm is $<\pi$,
 and beyond that it may or may not hold (compare, say, $\mathbb R$ and $\mathbb C$).
It is still easy to see that if the dimension $d$ of the unital algebra $\mathfrak A$ is $3$, then the algebra is either
 commutative or   isomorphic to the algebra of upper triangular matrices.

Now, upper triangular matrices already lead to some divergent Magnus series:

\begin{lemma}\plabel{lem:rys1}
Assume that $-\pi\leq\alpha\leq\pi$ and $\varepsilon\neq0$.
Let
\[M_1^{[\alpha,\varepsilon]}=\begin{bmatrix}
\frac{\pi-\alpha}2& -\frac{\pi+\alpha}2 \varepsilon\\&- \frac{\pi-\alpha}2
\end{bmatrix}
\qquad\text{and}\qquad
M_2^{[\alpha,\varepsilon]}=\begin{bmatrix}
\frac{\pi+\alpha}2 &
 \frac{\pi-\alpha}2 \varepsilon\\
&- \frac{\pi+\alpha}2
\end{bmatrix}.\]
Then the Magnus (BCH) series of the ordered measure $M_1^{[\alpha,\varepsilon]}\mathbf 1_{[0,1)}\boldsymbol. M_2^{[\alpha,\varepsilon]}\mathbf 1_{[1,2)}$ is not absolutely convergent.
\begin{proof}
If the series 
 $k\mapsto\mu_{k,\mathbb R}\left( M_1^{[\alpha,\varepsilon]}\mathbf 1_{[0,1)}\boldsymbol. M_2^{[\alpha,\varepsilon]}\mathbf 1_{[1,2)}\right)$
 were absolute convergent, then so would be
 $k\mapsto\mu_{k,\mathbb R}\left(\mathrm i\cdot(M_1^{[\alpha,\varepsilon]}\mathbf 1_{[0,1)}\boldsymbol. M_2^{[\alpha,\varepsilon]}\mathbf 1_{[1,2)})\right)$ in the complex(ified) algebra.
Assume that   $M^{[\alpha,\varepsilon]}$ is the absolute convergent sum of the Magnus series.
Then
\begin{multline*}\exp M^{[\alpha,\varepsilon]}=
\Rexp\left(\mathrm i\cdot(M_1^{[\alpha,\varepsilon]}\mathbf 1_{[0,1)}\boldsymbol. M_2^{[\alpha,\varepsilon]}\mathbf 1_{[1,2)})\right)=\\=\exp (\mathrm iM_1^{[\alpha,\varepsilon]})\exp (\mathrm iM_2^{[\alpha,\varepsilon]})=
\begin{cases}\begin{bmatrix}-1&-\frac{ \varepsilon(\pi^2+\alpha^2)(1+\mathrm e^{-\mathrm i \alpha }) }{\pi^2- \alpha^2}\\&-1\end{bmatrix}
&\text{for }  -\pi <  \alpha < \pi ,\\\\
\begin{bmatrix}-1&\pm\pi\mathrm i\varepsilon \\&-1\end{bmatrix}
&\text{for }   \alpha  =\pm \pi .
\end{cases}
\end{multline*}
In any case, the result is not diagonalizable, not even in the larger algebra of $2\times 2$ complex matrices.
It implies that $M^{[\alpha,\varepsilon]}$ must have Jordan form $\begin{bmatrix}a&1 \\&a\end{bmatrix}$ in the larger algebra of $2\times 2$ complex matrices
with $a\in2\pi\mathrm i\left(\mathbb Z+\frac12 \right)$.
This, however, contradicts to the fact that the terms of the Magnus (BCH) series 
 should be traceless in the larger algebra of $2\times 2$ complex matrices.
(The first term is known, and the other ones are sums of commutators.
Or, alternatively, this also follows directly from considering the homomorphisms given by the diagonal positions.)
\end{proof}
\end{lemma}

In the previous proof, instead of ``proper analysis'', only  linear algebra was used.
The following lemma, of more analytical nature,  indicates the oscillatory nature of the coefficients:

\begin{lemma}\plabel{lem:rys2}
In the previous lemma, for the terms of the Magnus expansion, $k\geq1$,
\begin{multline*}
\mu_{k,\mathrm R}\left( M_1^{[\alpha,\varepsilon]}\mathbf 1_{[0,1)}\boldsymbol. M_2^{[\alpha,\varepsilon]}\mathbf 1_{[1,2)}\right)=\begin{bmatrix}\pi& \\&-\pi\end{bmatrix}\cdot\delta_{k,1}+\\+
\begin{bmatrix}0&1 \\&0\end{bmatrix}\cdot
\begin{cases}(-1)^{\frac{k }2}\cdot-
\frac{2(\pi^2+ \alpha ^2)(1+\cos\alpha) \varepsilon}{\pi^2-\alpha ^2}+O(2^{-k})&\text{for  $k$ even, } -\pi <  \alpha < \pi ,\\\\
(-1)^{\frac{k+1}2}\cdot-
\frac{2(\pi^2+ \alpha^2)(\sin\alpha)\varepsilon}{
\pi^2- \alpha^2}+O(2^{-k})&\text{for $k$ odd, }  -\pi <  \alpha < \pi ,\\\\
O(2^{-k})&\text{for  $k$ even, }   \alpha  =\pm \pi ,\\\\
\pm(-1)^{\frac{k-1}2}2\pi\varepsilon+  O(2^{-k})&\text{for  $k$ odd, }    \alpha  =\pm \pi .
\end{cases}
\end{multline*}
In particular, the terms of the Magnus expansion limit to
$\begin{bmatrix}0& \pm c^{[\alpha]}_{\mathrm{odd}}\varepsilon\\&0\end{bmatrix}$
and
$\begin{bmatrix}0& \pm c^{[\alpha]}_{\mathrm{even}}\varepsilon\\&0\end{bmatrix}$
where $c^{[\alpha]}_{\mathrm{odd}},c^{[\alpha]}_{\mathrm{even}}\geq0 $, but at least one of them is nonzero.

\begin{proof}
Let
\[M^{[\alpha,\varepsilon]}(t)=\begin{cases}\begin{bmatrix}\pi t&
{\dfrac {t\varepsilon\,{\pi }  \left(
  (\pi^2+ \alpha^2) \left(
\cosh \left( t\pi\right) -\mathrm e^{  -  t  \alpha   }  \right)-
 2  \alpha  \pi\,\sinh \left( t\pi  \right)
 \right)
  }{ \left( {\pi }^{2}- \alpha^{2} \right)\sinh \left( t\pi  \right) }}
\\&-\pi t\end{bmatrix}
\\\hfill\text{for }-\pi <  \alpha < \pi ,\\
\begin{bmatrix}\pi t&\mp\pi t\varepsilon\left(\boldsymbol\beta(\pm2\pi t) \right) \\&-\pi t\end{bmatrix}
 \hfill\text{for }  \alpha =\pm \pi ;\\
\end{cases}\]
where we used the abbreviation $\boldsymbol\beta(x)=\frac{x}{\mathrm e^x-1}$.
Then $M^{[\alpha,\varepsilon]}(t)$ is meromorphic in $t$ with possible poles in $\mathrm i\mathbb Z$,
and it can be checked that
\[\exp M^{[\alpha,\varepsilon]}(t)=\Rexp\left(t\cdot(M_1^{[\alpha,\varepsilon]}\mathbf 1_{[0,1)}\boldsymbol. M_2^{[\alpha,\varepsilon]}\mathbf 1_{[1,2)})\right).\]
In fact, $M^{[\alpha,\varepsilon]}(t)$ extends to $t=0$ holomorphically with $M^{[\alpha,\varepsilon]}(0)=\Id_2$, and the exponential identity above also holds there.
As the continuation of $\log$ is unique near the identity,
\[M^{[\alpha,\varepsilon]}(t)=\log\Rexp\left(t\cdot(M_1^{[\alpha,\varepsilon]}\mathbf 1_{[0,1)}\boldsymbol. M_2^{[\alpha,\varepsilon]}\mathbf 1_{[1,2)})\right)
=\sum_{k=1}^\infty t^k\mu_{k,\mathrm R}\left( M_1^{[\alpha,\varepsilon]}\mathbf 1_{[0,1)}\boldsymbol. M_2^{[\alpha,\varepsilon]}\mathbf 1_{[1,2)} \right)\]
holds for $t\sim 0$.
Now, $\left(M^{[\alpha,\varepsilon]}(t)\right)_{1,2}$ has residue
\begin{align*}
&{\frac { ({\pi }^{2}+ \alpha ^2)  \varepsilon\, \left(
\mathrm i\right)
(1+\mathrm e^{-\mathrm i\alpha})  }{{
\pi }^{2}- \alpha ^{2}}}
&&\text{at } t=\mathrm i \text{ for }  -\pi <  \alpha < \pi ,\\
&
{\frac { ({\pi }^{2}+ \alpha ^2  )\varepsilon\, \left(-\mathrm  i\right)
(1+\mathrm e^{\mathrm i\alpha})  }{{
\pi }^{2}- \alpha ^{2}}}
&&\text{at } t=-\mathrm i \text{ for } -\pi <  \alpha < \pi ,\\
& \pm \pi \varepsilon
&&\text{at } t= \mathrm i \text{ for }   \alpha  =\pm \pi ,\\
& \pm \pi \varepsilon
&&\text{at } t= -\mathrm i \text{ for }    \alpha  =\pm \pi .
\end{align*}
Then the residues at $t=\pm\mathrm i$ give the macroscopic contributions in the power series,
the residues at $t=\pm 2\mathrm i$ give  $O(2^{-k})$, and the rest is altogether even smaller.
\end{proof}
\end{lemma}

\begin{example}\plabel{ex:rys}
In the (proof of) the previous lemma, the cases
$\alpha=0$ (the ``balanced cases'')
 and
$\alpha=\pm \pi $ (the ``totally unbalanced cases'')
are particularly instructive.

Here we find for $t\sim 0$, or even for $|t|<1$, or even just symbolically in $t$,
\[
\sum_{k=1}^\infty t^k\mu_{k,\mathrm R}\left( M_1^{[\alpha,\varepsilon]}\mathbf 1_{[0,1)}\boldsymbol. M_2^{[\alpha,\varepsilon]}\mathbf 1_{[1,2)}\right)
=\begin{cases}\begin{bmatrix}\pi t&
\varepsilon\pi t \tanh\left(\frac\pi2 t\right)
\\&-\pi t\end{bmatrix}
&\text{for } \alpha=0,\\\\
\begin{bmatrix}\pi t&\mp \pi t\varepsilon\left(\boldsymbol\beta(\pm2\pi t) \right) \\&-\pi t\end{bmatrix}
&\text{for }  \alpha =\pm \pi \\
\end{cases}
\]\[
=\begin{cases}\begin{bmatrix}\pi t&
\sum_{j=1}^\infty -4(-1)^j (1-2^{-2j})\zeta(2j)\varepsilon t^{2j}
\\&-\pi t\end{bmatrix}
&\text{for } \alpha=0,\\\\
\begin{bmatrix}\pi t&\mp \pi\varepsilon t + \pi^2 \varepsilon t^2\pm\sum_{j=1}^\infty (-1)^j2\pi\zeta(2j)\varepsilon t^{2j+1}  \\&-\pi t\end{bmatrix}
&\text{for }   \alpha  =\pm \pi .
\end{cases}
\]
(One may recall that as $1<2j\nearrow+\infty$, one has $\zeta(2j)\searrow1$, and even  $(1-2^{-2j})\zeta(2j)\searrow1$.)
\end{example}
The example(s) above are also instructive in illustrating that in the few cases when the
Magnus expansion of $\phi$ is explicitly known, it comes typically not from the direct evaluation of the
iterated integrals but from the knowledge of $\log(\Rexp(t\cdot\phi))$ for $t\sim 0$.
\begin{lemma}\plabel{lem:rys3}
Assume that the algebra of $2\times2$ (real or complex) upper triangular
matrices is endowed by an algebra norm $\|\cdot\|$ such that $\left\| \begin{bmatrix}1&\\&-1\end{bmatrix} \right\|=1$.
Then
\[\|M_1^{[\alpha,\varepsilon]}\|+\|M_2^{[\alpha,\varepsilon]}\|\leq \pi+\pi|\varepsilon|
 \left\| \begin{bmatrix}0&1\\&0\end{bmatrix}\right\|. \]
In the totally unbalanced cases,
\[\|M_1^{[-\pi,\varepsilon]}\|=\pi, \qquad \|M_2^{[-\pi,\varepsilon]}\|= \pi|\varepsilon|
 \left\| \begin{bmatrix}0&1\\&0\end{bmatrix}\right\|,\]
and
\[\|M_1^{[\pi,\varepsilon]}\|=\pi|\varepsilon|
 \left\| \begin{bmatrix}0&1\\&0\end{bmatrix}\right\|, \qquad \|M_2^{[\pi,\varepsilon]}\|= \pi.\]
\begin{proof}
These are trivial norm estimates.
\end{proof}
\end{lemma}
\snewpage

The main point is,
\begin{theorem}\plabel{thm:rys1}
Assume that the algebra of $2\times2$ (real or complex) upper triangular
 matrices is endowed by an $\|\cdot\|_{\ell^p}$ norm with $1\leq p\leq+\infty$,
 or even just with any algebra norm $\|\cdot\|$ such that $\left\| \begin{bmatrix}1&\\&-1\end{bmatrix} \right\|=1$.
Let $r>0$ and $\delta>0$. Then

(a) There exist an  example of type Lemma \ref{lem:rys1} such that
$\|M_1^{[\alpha,\varepsilon]}\|:\|M_2^{[\alpha,\varepsilon]}\|=r$ and
$\|M_1^{[\alpha,\varepsilon]}\|+\|M_2^{[\alpha,\varepsilon]}\|<\pi+\delta$.

(b) There exist a counterexample of BCH type $M_1\mathbf 1_{[0,1)}\boldsymbol. M_2\mathbf 1_{[1,2)}$ to
the convergence of the Magnus expansion such that
$\|M_1 \|:\|M_2 \|=r$ and
$\|M_1 \|+\|M_2 \|=\pi+\delta$, and the coefficients of the Magnus expansion are unbounded.
\begin{proof}
(a) We can choose a small $\varepsilon>0$ such that
$r\in\left[ \varepsilon \left\| \begin{bmatrix}0&1\\&0\end{bmatrix}\right\|,
\varepsilon^{-1} \left\| \begin{bmatrix}0&1\\&0\end{bmatrix}\right\|^{-1}
\right]$
and $\pi \varepsilon\left\| \begin{bmatrix}0&1\\&0\end{bmatrix}\right\| <\delta$.
According the previous lemma, the first relation implies that
as we pass from one totally unbalanced case to the other, the norm ratio $r$ will be realized.
The second relation implies that the cumulative norm is still small.

(b) We can linearly upscale (a).
\end{proof}
\end{theorem}
\begin{remark}\plabel{rem:rys}
The $\|\cdot\|_{\ell^p}$ norms with $p=1,2,+\infty$ can be computed relatively explicitly.
For example, regarding the cumulative norm
\[\|M_1^{[\alpha,\varepsilon]}\|_{\ell^p}+\|M_2^{[\alpha,\varepsilon]}\|_{\ell^p}=\pi+\pi|\varepsilon|\qquad\text{if }p=1,\infty;\]
and
\begin{multline*}
\|M_1^{[\alpha,\varepsilon]}\|_{\ell^2}+\|M_2^{[\alpha,\varepsilon]}\|_{\ell^2}=\\
=\frac12\pi|\varepsilon|+ \sqrt{\left(\frac{\pi-\alpha}2\right)^2+\left(\frac{\pi+\alpha}2\right)^2\frac{|\varepsilon|^2}4}
+ \sqrt{\left(\frac{\pi+\alpha}2\right)^2+\left(\frac{\pi-\alpha}2\right)^2\frac{|\varepsilon|^2}4}
\\
\leq \pi+\pi|\varepsilon|
.\eqedremark
\end{multline*}
\end{remark}
\begin{remark}\plabel{rem:rys11}
The counterexamples in the previous theorem fit to the line 
of Wei \cite{W} in using the argument of parabolicity (but we use it in a slightly more sophisticated way),
of Michel \cite{Mi} in finding a whole range of balanced and unbalanced counterexamples (but we consider a different setting),
and of Vinokurov \cite{V} in finding relatively simple counterexamples (but we use an even simpler algebra and
 exhibit a bigger range of  counterexamples).
\qedremark
\end{remark}
\snewpage
I prefer to think about the construction of Lemma \ref{lem:rys1} as the ``Minimal Examples''.
This is not about any mathematical minimality property but that the construction realizes counterexamples in a quite minimal setting.
There are variants of it, but one cannot go below cumulative norm $\pi$:
\begin{theorem}\plabel{thm:rys2}
Assume that $\mathfrak A$ is the algebra of the $2\times2$ upper triangular matrices
 (real or complex), endowed by a Banach algebra norm $|\cdot|_{\mathfrak A}$.
Assume that $\phi$ is a continuous $\mathfrak A$-valued measure of cumulative norm $\leq\pi$.
Then the Magnus series of $\phi$ is (absolute) convergent.
\begin{proof}
Using standard Neumann series arguments, it is easy to see that the absolute value of the diagonal elements is dominated by the norm.
More precisely, if $X\in \mathfrak A$, then $|(X)_{1,1}|\leq |X|_{\mathfrak A}$ and $|(X)_{2,2}|\leq |X|_{\mathfrak A}$.
If $\int|\phi|_{\mathfrak A}<\pi$, then in the complex(ified) setting for $t\in\Dbar(0,1+\varepsilon)$,
 $M(t)=\log\Rexp(t\cdot\phi)$ exists as one can apply the definition
 $\log A=\int_{\lambda\in{[0,1]}} \frac{A-1}{\lambda+(1-\lambda)A}\,\mathrm d\lambda$
 based on the critical behavior in the diagonals.
The power series expansion of $M(t)$ around $t=0$ gives the convergent Magnus expansion for $t=1$.
In fact, the same applies even if $\int|\phi|_{\mathfrak A}=\pi$ but $\left|\int(\phi)_{1,1}\right|<\pi$
 and $\left|\int(\phi)_{2,2}\right|<\pi$.
Let us now consider the measure $\tilde \phi$ such that $\phi=\frac12(\int(\phi)_{1,1}+\int(\phi)_{2,2})\cdot\Id_2+\tilde \phi$.
I. e. $\tilde\phi$ is the traceless modification of $\phi$.
The  Magnus expansion of $\tilde \phi$ is equiconvergent to the  Magnus expansion of $ \phi$.
This traceless modification is convergent as the integrated modified diagonals get off from the boundary of $\Dbar(0,\pi)$
except if  $C=\int(\phi)_{1,1}=-\int(\phi)_{2,2}$ with $|C|=\int|\phi|_{\mathfrak A}=\pi$.
It is sufficient to consider this latter case.
Reparametrizing the measure $\phi$ in terms of variation, and multiplying it a complex unit vector, we can assume that
$\phi$ is supported on $[0,\pi)$ and
$\phi=\begin{bmatrix}\mathbf 1_{[0,\pi)}&\mathbf c\\&-\mathbf 1_{[0,\pi)}\end{bmatrix}$ and $|\psi|_{\mathfrak A}=\mathbf 1_{[0,\pi)}$.
Now, $\mathbf c$ may be essentially constant i. e. a scalar times   $\mathbf 1_{[0,\pi)}$.
In this case the Magnus expansion is convergent.
If this is not so, then there are elements
$\begin{bmatrix}1&a\\&-1\end{bmatrix}$  and $\begin{bmatrix}1&b\\&-1\end{bmatrix}$ of unit norm
such that $a\neq b$.
But then the norm of $n\begin{bmatrix}0& (b-a)\\&0\end{bmatrix}
=\left(\begin{bmatrix}1&a\\&-1\end{bmatrix}\begin{bmatrix}1&b\\&-1\end{bmatrix}\right)^n-\begin{bmatrix}1&\\&1\end{bmatrix}$
would be bounded by $2$, which is an absurdum.
\end{proof}
\end{theorem}
If dimension $4$ is allowed for an algebra $\mathfrak A$,
then there are already counterexamples to cumulative $\pi$ with respect to the convergence of the Magnus expansion.
Indeed, in the case of real $2\times 2$ matrices with the $\ell^2$ operator norm
the Moan--Schäffer example (see Schäffer \cite{Scha}, Moan \cite{Ma}, Moan, Niesen \cite{MN})
 or the Magnus critical example (see Part II \cite{L2} for a detailed discussion) are like that.

Taking $\oplus\mathbb R$ and the joint maximum norm, we can always increase the dimension of the
Banach algebras while leaving the cumulative norm of the counterexample the same.
However, we have promised a series of counterexamples with cumulative norm $\searrow2$.
\snewpage
~

\subsection{Higher dimensional counterexamples}
\plabel{ss:high}
~\\

Let us start with a particular example, which we will generalize.
\begin{example}
\plabel{ex:div}
Here we will construct an example for the divergence of Magnus expansion
 which is of mBCH type, using $5\times5$ real matrices of cumulative norm $5\cdot\frac2{5-1}=\frac52$.
We will use the matrices
\[M^{(5)}_1=\frac24\left[\begin{matrix}0&1&1&1&1\\&0&&&\\&&0&&\\&&&0&\\&&&&0\end{matrix}\right],\,\,\,
M^{(5)}_2=\frac24\left[\begin{matrix}0&&&&\\-1&0&1&1&1\\&&0&&\\&&&0&\\&&&&0\end{matrix}\right],\,\,\,
M^{(5)}_3=\frac24\left[\begin{matrix}0&&&&\\&0&&&\\-1&-1&0&1&1\\&&&0&\\&&&&0\end{matrix}\right],\]
\[M^{(5)}_4=\frac24\left[\begin{matrix}0&&&&\\&0&&&\\&&0&&\\-1&-1&-1&0&1\\&&&&0\end{matrix}\right],\,\,\,
M^{(5)}_5=\frac24\left[\begin{matrix}0&&&&\\&0&&&\\&&0&&\\&&&0&\\-1&-1&-1&-1&0\end{matrix}\right].\]

The measure we will consider is
$\psi_5=M^{(5)}_1\mathbf 1_{[0,1)}\boldsymbol.
M^{(5)}_2\mathbf 1_{[1,2)}\boldsymbol.
M^{(5)}_3\mathbf 1_{[2,3)}\boldsymbol.
M^{(5)}_4\mathbf 1_{[3,4)}\boldsymbol.
M^{(5)}_5\mathbf 1_{[4,5)}$.
For the matrices we will use the $\ell^1$ operator norm $\|\cdot\|_{\ell^1}$.
Then the cumulative norm of $\psi_5$ is
\[\int\|\psi_5\|_{\ell^1}=\|M^{(5)}_1\|_{\ell^1}+\ldots+\|M^{(5)}_5\|_{\ell^1}=5\cdot\frac24=\frac52.\]
On the other hand,
\begin{align*}
\Rexp(\psi_5)&=(\exp M^{(5)}_1)\ldots(\exp M^{(5)}_5)
=(\Id_5+M^{(5)}_1)\ldots(\Id_5+M^{(5)}_5)
\\
&= \left[ \begin {matrix} -{\dfrac {33}{32}}&-{\dfrac {41}{32}}&-{
\dfrac {21}{32}}&{\dfrac {9}{32}}&{\dfrac {27}{16}}\\\\  -
{\dfrac {27}{16}}&-\dfrac{3}{16}&-{\dfrac {7}{16}}&\dfrac3{16}&{\dfrac {9}{8}}
\\\\  -{\dfrac {9}{8}}&-{\dfrac {9}{8}}&\dfrac38&\dfrac18&\dfrac34
\\\\  -\dfrac34&-\dfrac34&-\dfrac34&\dfrac34&\dfrac12\\\\
 -\dfrac12&-
\dfrac12&-\dfrac12&-\dfrac12&1\end {matrix} \right]
=
 \left[ \begin {matrix} 1&-7&-\sqrt {15}&5&-2
\\  -1&-2&2\,\sqrt {15}&5&-1\\  1&8&0
&5&0\\  -1&-2&-2\,\sqrt {15}&5&1\\
1&-7&\sqrt {15}&5&2\end {matrix} \right]
\cdot\\&\qquad\cdot
\underbrace{\left[ \begin {matrix} 1&0&0&0&0\\  0&{\frac {
61}{64}}&-{\frac {5\,\sqrt {15}}{64}}&0&0\\  0&{
\frac {5\,\sqrt {15}}{64}}&{\frac {61}{64}}&0&0\\  0&0
&0&-1&1\\  0&0&0&0&-1\end {matrix} \right]}_{=F_5}
 \cdot
 \left[ \begin {matrix} 1&-7&-\sqrt {15}&5&-2
\\  -1&-2&2\,\sqrt {15}&5&-1\\  1&8&0
&5&0\\  -1&-2&-2\,\sqrt {15}&5&1\\
1&-7&\sqrt {15}&5&2\end {matrix} \right]^{-1}
\end{align*}
informs us about the real Jordan form $F_5$ of $\Rexp(\psi_5)$.
We see that the geometric multiplicity of the eigenvalue $-1$ is $1$,
therefore it cannot be the exponential of a real matrix $M$.
(The eigenvalue $-1$ is not forbidden for the exponential of a real matrix but
its geometric multiplicity should be even.)
On the other hand, the sum $M$ of the Magnus series should be exactly a matrix like that.
\qedexer
\end{example}
That will be our general strategy:
 we will consider measures of real matrices with time-ordered exponential with
 eigenvalue $-1$ of geometric multiplicity $1$.
In order to establish this behaviour we, in fact, do not have to deal with
 the Jordan form, a simple rank computation for $\Rexp(\phi)+\Id$ is sufficient.

 Now, we start the general construction.
Let $Q^{(n)}_s$ be the $n\times n$ real matrix such that its elements ($i$th row, $j$th column)
are given by
\[\left(Q^{(n)}_u\right)_{i,j}=\delta_{i,u}\sgn(j-u)\]
(where the Kronecker's delta notation is used).
\begin{lemma}\plabel{lem:parprod}
For any $1\leq k\leq n$, and any scalar $s$,
\begin{multline*}(\exp sQ^{(n)}_1)\ldots(\exp sQ^{(n)}_k)=(\Id_n+ sQ^{(n)}_1)\ldots(\Id_n+ sQ^{(n)}_k)
=\\=
\left[\begin{array}{c|c}
A^{(k)} (s)+B^{(k)}(s)&C^{(k,n-k)}(s)\\\hline 0_{(n-k)\times k}&\Id_{n-k}
\end{array}\right];
\end{multline*}
such that
$A^{(k)}(s)$ is a $k\times k$ matrix whose elements ($i$th row, $j$th column)
are given by
\[\left(A^{(k)}(s)\right)_{i,j}=
\begin{cases}
(s+1)^{j-i+1}-(s+1)^{j-i-1}&\text{if }\quad i<j,\\
s+1&\text{if }\quad i=j,\\
0&\text{if }\quad i>j;
\end{cases}\]
$B^{(k)}(s)$ is a $k\times k$ matrix whose elements ($i$th row, $j$th column)
are given by
\[\left(B^{(k)}(s)\right)_{i,j}=(s+1)^{k-i}-(s+1)^{k-i+1};\]
$C^{(k,n-k)}(s)$ is a $k\times (n-k)$ matrix whose elements ($i$th row, $j$th column)
are given by
\[\left(C^{(k,n-k)}(s)\right)_{i,j}=-(s+1)^{k-i}+(s+1)^{k-i+1}.\]
(The elements $B^{(k)}(s)$  and $C^{(k,n-k)}(s)$ depend only on the row number.)
\begin{proof}
One can prove this by induction in $k$.
\end{proof}
\end{lemma}
\begin{lemma}\plabel{lem:parprod2}
(a) For any $1\leq n$, and any scalar $s$,
\[
P^{(n)}(s):=(\exp sQ^{(n)}_1)\ldots(\exp sQ^{(n)}_n)=(\Id_n+ sQ^{(n)}_1)\ldots(\Id_n+ sQ^{(n)}_n)
=A^{(n)} (s)+B^{(n)}(s)
\plabel{eq:proder}
\]

(b) If $s+1\neq \pm 1$, then he eigenvalues of product matrix $P^{(n)}(s)$ have geometric multiplicity at most $1$.
If $s=0$, then $\lambda=s+1=1$ has geometric multiplicity $n$;
if $s=2$, then $\lambda=s+1=-1$ has geometric multiplicity $n-1$.

(c) If $s=2+\frac2{n-1}$, then $-1$ is an eigenvalue of $P^{(n)}(s)$.
\begin{proof}
(a) This is an immediate consequence of the previous lemma.

(b) Let us consider the matrix $P^{(n)}(s)-\lambda\Id_n=(A^{(n)} (s)-\lambda\Id_n) +B^{(n)}(s)$.

If $\lambda\neq s+1$, then $(A^{(n)} (s)-\lambda\Id_n)$ is an invertible triangular matrix
of full rank, while $B^{(n)}(s)$ is rank $0$ or $1$. Therefore the rank of the sum is $n-1$ of $n$.

If $\lambda= s+1\neq -1$, then (we also know that $s+1\neq1$) we have
$(s+1)^2-1\neq 0$, thus
the column space $(A^{(n)} (s)-\lambda\Id_n)$ contains the
column vectors whose last entry is $0$,
and $(A^{(n)} (s)-\lambda\Id_n)$ also has a $0$ column.
Meanwhile the columns of $B^{(n)}(s)$
are uniformly a column vector, whose last entry is $1-(s+1)=-s\neq0$.
Then it is easy to see that the columns of $(A^{(n)} (s)-\lambda\Id_n) +B^{(n)}(s)$
must be independent.

The special cases $\lambda= s+1=\pm1$ can be checked directly.

(c) using (a) one can check that the column vector with uniform entries $1$ is an eigenvector with eigenvalue $-1$.
\end{proof}
\end{lemma}
Now it is clear what to do, set $M^{(n)}_k=\frac2{n-1}Q^{(n)}_k$, and
\[\psi_n=M^{(n)}_1\mathbf 1_{[0,1)}\boldsymbol.
\ldots\boldsymbol.
M^{(n)}_k\mathbf 1_{[k-1,k)}\boldsymbol.
\ldots\boldsymbol.
M^{(n)}_n\mathbf 1_{[n-1,n)}.\]
Then we obtain
\begin{theorem}\plabel{th:parprod} Let $n\geq2$.
The cumulative norm of $\psi_n$ in the $\ell^1$ operator norm $\|\cdot\|_{\ell^1}$ of $n\times n$ matrices is
\[\int \|\psi_n\|_{\ell^1}=\frac{2n}{n-1},\]
while the Magnus expansion of $\psi_n$ is not absolutely convergent.
\begin{proof}
The statements about the cumulative norms are straightforward.
As $-1$ is an eigenvalue of geometric multiplicity $1$ of $\Rexp (\psi_n)$,
it cannot an exponential of a real matrix.
Meanwhile, the sum $M$ of the Magnus series should be  a matrix like that.
\end{proof}
\end{theorem}

\begin{remark}\plabel{rem:parprod}
One can see that the eigenvalues of $\Rexp (\psi_n)$ are on the unit circle,
with the only nontrivial Jordan block
$\begin{bmatrix}-1&1\\&-1\end{bmatrix}$ in complex form
(corresponding to column vectors of arithmetic progressions).

Indeed, one can argue as follows.
Let $\mathbf v$ be the column vector of length $n$ containing only entries $1$.
Then the matrix $\Id_n-\frac1n\mathbf v\mathbf v^\top$ defines a positive semidefinite
 quadratic form $S_n$ on $\mathbb R^n$, but it descends to a positive definite
 quadratic form $\tilde S_n$ on the factor space $\mathbb R^n/\mathbb R\mathbf v=\tilde{\mathbb R}^n$.
Then $P_n=\Rexp(\psi_n)$ leaves $S_n$ invariant, and it also descends to linear map $\tilde P_n $ on $\tilde{\mathbb R}^n$,
 which is therefore orthogonal with respect to $\tilde S_n$.
Thus the eigenvalues of $\tilde P_n$ are unit complex numbers with trivial Jordan blocks.
Getting back to $\mathbb R^n$, can gives one extra eigenvalue from $\mathbb R\mathbf v$
which we know to be $-1$, and it can give only one extra nontrivial Jordan block, which
must be as indicated because the effect of $P_n$ is easy to check on column vectors of arithmetic progressions.
\qedremark
\end{remark}
\snewpage
\begin{commenty}
\begin{example}\plabel{ex:parprod}
\[\Rexp(\psi_2)=\left[ \begin {array}{cc} -3&2\\ \noalign{\medskip}-2&1\end {array}\right]=
\left[ \begin {array}{cc} 2&-\frac12\\ \noalign{\medskip}2&\frac12\end {array} \right]
\left[ \begin {array}{cc} -1&1\\ \noalign{\medskip}0&-1\end {array}\right]
\left[ \begin {array}{cc} 2&-\frac12\\ \noalign{\medskip}2&\frac12\end {array} \right]^{-1}
;\]
\[
\Rexp(\psi_3)=
\left[ \begin {array}{ccc} -2&-1&2\\ \noalign{\medskip}-2&0&1
\\ \noalign{\medskip}-1&-1&1\end {array} \right]
=
\left[ \begin {array}{ccc} 1&3&-1\\ \noalign{\medskip}-1&3&0
\\ \noalign{\medskip}1&3&1\end {array} \right]
\left[ \begin {array}{ccc} 1&0&0\\ \noalign{\medskip}0&-1&1
\\ \noalign{\medskip}0&0&-1\end {array} \right]
\left[ \begin {array}{ccc} 1&3&-1\\ \noalign{\medskip}-1&3&0
\\ \noalign{\medskip}1&3&1\end {array} \right]^{-1}
;\]
\begin{multline*}
\Rexp(\psi_4)=
 \left[ \begin {array}{cccc} -{\frac {115}{81}}&-{\frac {106}{81}}&-{
\frac {10}{81}}&{\frac {50}{27}}\\ \noalign{\medskip}-{\frac {50}{27}}
&-{\frac {5}{27}}&-{\frac {2}{27}}&{\frac {10}{9}}
\\ \noalign{\medskip}-{\frac {10}{9}}&-{\frac {10}{9}}&\frac59&\frac23
\\ \noalign{\medskip}-\frac23&-\frac23&-\frac23&1\end {array} \right]
= \left[ \begin {array}{cccc} 7&\sqrt {5}&4&-\frac32\\ \noalign{\medskip}-3
&-3\,\sqrt {5}&4&-\frac12\\
\noalign{\medskip}-3&3\,\sqrt {5}&4&\frac12
\\ \noalign{\medskip}7&-\sqrt {5}&4&\frac32\end {array} \right]
\cdot\\\cdot
\left[ \begin {array}{cccc} {\frac {79}{81}}&-{\frac {8\,\sqrt {5}}{
81}}&0&0\\ \noalign{\medskip}{\frac {8\,\sqrt {5}}{81}}&{\frac {79}{81
}}&0&0\\ \noalign{\medskip}0&0&-1&1\\ \noalign{\medskip}0&0&0&-1
\end {array} \right]
\left[ \begin {array}{cccc} 7&\sqrt {5}&4&-\frac32\\ \noalign{\medskip}-3
&-3\,\sqrt {5}&4&-\frac12\\
\noalign{\medskip}-3&3\,\sqrt {5}&4&\frac12
\\ \noalign{\medskip}7&-\sqrt {5}&4&\frac32\end {array} \right]^{-1}.
\eqedexer
\end{multline*}
\end{example}
\end{commenty}

For $d\geq3$, let $\mathrm C^{\{\{d\}\}}_{\mathbb K}$ denote the infimum
 of the cumulative norms of $\mathfrak A$ valued ordered measures whose Magnus expansion is not absolutely convergent
 and $\mathfrak A$ is a $d$ dimensional Banach algebra over $\mathbb K$.
We know that $\mathrm C^{\{\{3\}\}}_{\mathbb K}=\pi$, and $\mathrm C^{\{\{d\}\}}_{\mathbb K}$
 is (possibly not strictly) decreasing in $d$, but $\geq2$.
By Theorem \ref{th:parprod}, for dimension $d=n^2$ ($n\geq2$), the continuous measure $\psi_n$ provides a counterexample
 to the convergence of the Magnus expansion.
Thus, for $d\geq4$,
\begin{equation}
\mathrm C^{\{\{d\}\}}_{\mathbb R}\leq 2+\frac2{ \lfloor\sqrt d \rfloor-1}.
\plabel{eq:ram}
\end{equation}
(Meanwhile, $\mathrm C^{\{\{d\}\}}_{\mathbb R}\geq \mathrm C^{\{\{d\}\}}_{\mathbb C}\geq \mathrm C^{\{\{2d\}\}}_{\mathbb R}$ is trivial.)

Note, however, that the cumulative of $\psi_2$ is $4$ and the   cumulative of $\psi_3$ is $3$.
So, these counterexamples get lower cumulative norms than $\pi$ only for dimensions $d\geq9$.

For    $4\leq\dim_{\mathbb K}\mathfrak A\leq8$,
 it would be interesting to see counterexamples with cumulative norm less than $\pi$, if they exist.
In particular, most prominently, it is not clear how far one can go when $\mathfrak A$ is isomorphic
 to the algebra of real $2\times2$ matrices (allowing other norms than the $\ell^2$ operator norm).
\snewpage
~

\subsection{Relationship to the free $L^1$ mBCH counterexamples}
\plabel{ss:rel}
~\\

In Part I \cite{L1} we have considered counterexamples of mBCH type for the convergence of the  Magnus expansion.
There, the tautological measure $2\cdot\mathrm Z^1_{[0,1)}$ was replaced by
\begin{multline*}
\hat\psi_n=\log\Rexp(2\cdot\mathrm Z^1_{[0,\frac1n)} )\mathbf1_{[0,1)}\boldsymbol.\ldots\\
\ldots\boldsymbol.
\log\Rexp(2\cdot\mathrm Z^1_{[\frac{k-1} n,\frac{k}n )})\mathbf1_{[k-1,k)}\boldsymbol.\ldots\boldsymbol.
\log\Rexp(2\cdot\mathrm Z^1_{[\frac{n-1} n,\frac{n}n)} )\mathbf1_{[n-1,n)}
\end{multline*}
(with $n\geq2$).
In that the cumulative norm of the measure increases only slightly, but the time-ordered exponentals
 will be the same in $\mathrm F^{1,\mathrm{loc}}([0,1))$ but not lying in $\mathrm F^{1 }([0,1))$.
These counterexample can be understood in the algebra finitely generated by the
 $\log\Rexp(2\cdot\mathrm Z^1_{[\frac{k-1} n,\frac{k}n )})$
 but which is, nevertheless infinite dimensional and not particularly manageable.
Also, the cumulative norm
 $\int|\hat\psi_n|=n\Theta\left(\frac2n\right)=2+\frac2n+O\left(\frac1{n^2}\right) $
 and the cumulative norm $\int\|\psi_n\|_{\ell^1}=2+\frac2{n-1}$ are quite comparable.
According to these, the counterexamples $\psi_n$ compare preferably to the $\hat\psi_n$.

Here we would like to argue that the $\psi_n$ and $\hat\psi_n$ closely related to each other.
The point is that  $\psi_n$ can be obtained from $\hat\psi_n$ by ``reducing'' the algebra $\mathrm F^{1 }([0,1))$.
Although the following discussion can be made completely precise, for the sake of ease, we keep it informal.
In particular, we will sometimes pretend that $\mathrm F^{1 }([0,1))$ is the superposition of
 elements $X_{t_1}\ldots X_{t_k}$ with $t_i\in[0,1)$, even if this is not so simple.
 The precise idea is that $\mathrm F^{1 }([0,1))$ can be subjected to contractive homomorphisms but which
 can be set up so that the terms of the Magnus series are kept seen relatively large.

The first thing one can do is to replace $X_{t_1}\ldots X_{t_k}$ by
 $\mathsf{a}^{\asc(t_1,\ldots,t_k)}\mathsf{d}^{\des(t_1,\ldots,t_k)}X_{t_1,t_2}$,
 i.~e. the internal structure of the expression gets ignored, only the number of ascents and descents get recorded.
Actually, it is sufficient to have $\mathsf{a}\equiv1$ and $\mathsf{d}\equiv-1$ here; this setup still keeps
 the size of the terms of the Magnus expansion.
This leads to an ``abstract composition kernel approach''.
What can make this more down-to-earth is to consider it as a representation.
We have the representation space generated by superpositions of $Y_t$ ($t\in [0,1)$),
where the representation rule is $X_{t_1} Y_{t_2}=Y_{t_1}$ for $t_1<t_2$ and $X_{t_1} Y_{t_2}=-Y_{t_1}$ for $t_1>t_2$.
If the superpositions of $Y_t$ is coded as an  $L^1$ function $f(t)$, then the effect of $Z_{[a,b)}$ on  $f(t)$ leads to $\tilde f(t)$
where
\[\tilde f(t)=\int_{s\in[0,1)} \underbrace{\chi_{[a,b)}(t) \sgn(s-t)}_{=K_{[a,b)}(s,t)}\,f(s)\,\mathrm ds.\]
Thus $Z_{[a,b)}$ gets represented by kernel $K_{[a,b)}(s,t)$ ($L^1$ in $s$, $L^\infty$ in $t$).
This representation is still quite strange, and indeed, there is no loss in norm with respect to terms of the Magnus expansion.
\snewpage

What one can make the situation tamer is the introduction of additional eliminating relations $X_{t_1}X_{t_2}=0$ for $t_1,t_2\in[\frac{k-1}n,\frac{k}n )$.
This \textit{will} improve the convergence of the Magnus expansion but hopefully not so much.
This is compatible to the representation process of the previous paragraph.
It leads to the kernels $K^{[1/n]}_{[a,b)}(s,t)$ where the values on the squares $[\frac{k-1}n,\frac{k}n )\times [\frac{k-1}n,\frac{k}n ) $
are killed.
Then, the image of the time-ordered exponential of $2\cdot\mathrm Z^1_{[0,1)}$, and of $\hat \psi_n$ and the computations 
with the latter one can be
viewed in terms of kernels which are linear combinations of characteristic functions of $[\frac{k-1}n,\frac{k}n )\times [\frac{l-1}n,\frac{l}n )$.
I. e. in terms of $n\times n$ matrices.
Indeed, this is exactly the process how $\psi_n$ was obtained from $\hat \psi_n$.
\snewpage

\section{Convergence for finite dimensional Banach algebras}
\plabel{sec:ConFin}

We recall certain estimates due to Rogers, which are standard material.
Let us consider any centrally symmetric compact convex body $H$ in the $n$-dimensional space $\mathbb R^n$.
Let $\vartheta(H)$ be the infimum of the covering density of $\mathbb R^n$ by translates
 of $H$; and let $\vartheta_L(H)$ be the same with respect to lattice coverings.
Then, by Rogers \cite{R}, for $n\geq3$,
\begin{align}
\vartheta(H)\leq \vartheta_L(H)\leq&\min_{0<\eta<1/n}(1+\eta)^n(1+n\log(1/\eta))\plabel{eq:r1}\\
&=-n{  W}_{-1} \left(-{\tfrac1n}\right) \left( 1-{\frac {1}{n {  W}_{-1} \left(-{\tfrac1n}\right)} } \right) ^{n+1}\notag\\
<& \left(1+\frac1{n\log n}\right)^n(1+n\log(n\log n))\plabel{eq:r2}\\
<& n\log n +n\log\log n+2n+1 \plabel{eq:r3}\\
<& n\log n +n\log\log n+5n.\plabel{eq:r4}
\end{align}
Here  \eqref{eq:r1} is the proper result of \cite{R} (which is expressed explicitly using the Lambert $W_{-1}$ function);
 \eqref{eq:r2} reflects the choice $\eta=\frac1{n\log n}$; and
  while, for example,  \eqref{eq:r3} is still true,
 \eqref{eq:r4}  is the much quoted estimate;
 see G. Fejes Tóth \cite{F} for more on this.
The main point is that there is a dimension-dependent but otherwise universal quantity $\vartheta_n$
 (which can be chosen as any expression on the RHS),
 which is nearly linear in $n$, estimating the minimal translative covering densities from above.
Moreover, as it is explained in Rogers, Zong \cite{RZ}, if $0<r<1$, then $H$ can be covered by at
 most $(1+r^{-1})^n\vartheta(H)\leq (1+r^{-1})^n\vartheta_n$
 translates of $rH$ (i. e. by homothetical copies of $H$ with ratio $r$).
For a recent review on these topics, see Naszódi \cite{N}.
\begin{theorem}\plabel{thm:gain2}
Assume that $\mathfrak A$ is Banach algebra of finite real dimension $n$.
Assume that $\phi$ is an ordered measure of cumulative norm $\int|\phi|=\omega$.
Let $\lambda\in[0,1]$.
Then
\begin{multline}
\left|\mu_{2,\mathrm R}(\phi)\right|=
\left|\int_{t_1,t_2\in[0,1]} \lambda^{\asc(t_1,t_2)} (\lambda-1)^{\des(t_1,t_2)}\phi(t_1)\phi(t_2)\right|
\\
\leq\frac{\omega^2}2\left(1- {\frac {{2}^{1-n}}{n} \left( {\frac {1-\frac2n}{\mathrm e}} \right)  }\frac1{\vartheta_n}
\min(\lambda,1-\lambda)\right).
\plabel{eq:gar3}
\end{multline}

\begin{proof}
The measure $\phi$ can be approximated by measures of mBCH type.
Therefore it is sufficient to prove the statement in the case
\[\phi=u_1\mathbf1_{ [\tau_0,\tau_1)}\boldsymbol.\ldots \boldsymbol.u_k\mathbf1_{\mathbf [\tau_{k-1},\tau_k)}\]
 such that $u_i\in\mathfrak A$, $|u_i|=1$, $0=\tau_0<\tau_1<\ldots<\tau_k=\omega$.
We also write $\phi(t)=u(t)\mathbf 1_{[0,\omega)}$, where $u(t)$ is the corresponding piecewise constant  function.
\snewpage

Let us apply the argument of  Rogers, Zong \cite{RZ} where $ H$ is the closed unit ball $K$ of $\mathfrak A\simeq\mathbb R^n$.
Then $K$ can be covered by at most $s= (1+r^{-1})^n\vartheta_n $ many copies of $rK$, say $K_1,\ldots, K_{\lfloor s\rfloor}$.
Let $L_1=K\cap K_1$, and $L_i=(K\cap K_i )\setminus  (K_1\cup \ldots \cup K_{i-1})$.
Let $N_i$ be the union of those $[\tau_{p-1},\tau_p)$ such that $u_p\in L_i$.
(I. e., formally, $N_i=u^{-1}(L_i)$.)
Let $v_i$ be the Lebesgue measure of $N_i$.
Clearly, $v_1+\ldots+v_s=\omega$.
We cut $N_i$ into a lower part $N_i^-$ and lower part $N_i^+$, each of them of Lebesgue measure $v_i/2$.
Let $P_i^-:[0,v_i/2)\rightarrow N_i^-$ and $P_i^+:[0,v_i/2)\rightarrow N_i^+$
 be the corresponding monotone increasing measure preserving transformations.
It is easy to see that $P_i^+(q)-P_i^-(q)\geq v_i/2$ for any $q$.
On the other hand, $|u(P_i^+(q))-u(P_i^-(q))|\leq 2r$, by the construction of $L_i\subset K_i$.

One can write
\begin{multline*}
\mu_{2,\mathrm R}(\phi)
=\int_{t_1,t_2\in[0,\omega]} \lambda^{\asc(t_1,t_2)} (\lambda-1)^{\des(t_1,t_2)} u(t_1)u(t_2)\, \mathrm dt_1\,\mathrm dt_2\\
=\sum_i\int_{q\in [0,v_i/2),t_2\in[0,\omega]} \lambda^{\asc(P_i^-(q),t_2)} (\lambda-1)^{\des( P_i^-(q),t_2)} u(P_i^-(q))u(t_2)
+\\+ \lambda^{\asc(P_i^+(q),t_2)} (\lambda-1)^{\des( P_i^+(q),t_2)} u(P_i^-(q))u(t_2)
\, \mathrm dq\,\mathrm dt_2.
\end{multline*}

The norm of $(\lambda^{\asc(t_{11},t_2)} (\lambda-1)^{\des(t_{11},t_2)} u(t_{11})+
\lambda^{\asc(t_{12},t_2)} (\lambda-1)^{\des(t_{12},t_2)} u(t_{12})
)u(t_2)$ can be estimated
$(\lambda^{\asc(t_{11},t_2)} (1-\lambda)^{\des(t_{11},t_2)}\cdot1  +
\lambda^{\asc(t_{12},t_2)} (1-\lambda)^{\des(t_{12},t_2)}\cdot1  )\cdot1
) $, but this leads only to the trivial estimate
\[|\mu_{2,\mathrm R}(\phi)|\leq\frac12\omega^2.\]
However, when $t_{11}<t_2<t_{12}$ or $t_{12}<t_2<t_{11}$ holds, then there is gain,
 as a component $\min(\lambda,1-\lambda)\cdot 2$ in the estimate can be replaced by
 $\min(\lambda,1-\lambda)\cdot | u(t_{11})-u(t_{12})|$.
This is a gain $\min(\lambda,1-\lambda)\cdot (2-| u(t_{11})-u(t_{12})|)$ in the estimate.
Therefore, we have a better overall estimate
\begin{equation}
|\mu_{2,\mathrm R}(\phi)|\leq\frac{\omega^2}2- \sum_i\min(\lambda,1-\lambda)\cdot(2-2r)\cdot\left(\frac{v_i}2\right)^2.
\plabel{eq:gar1}
\end{equation}
The inequality between the arithmetic and square means implies that
$\sum (v_i)^2\geq \frac{\omega^2}{\lfloor s \rfloor}\geq \frac{\omega^2}s$.
Therefore \eqref{eq:gar1} implies
\begin{equation}
|\mu_{2,\mathrm R}(\phi)|\leq
\frac{\omega^2}2\left(1- \frac{1-r}s\min(\lambda,1-\lambda)\right).
\plabel{eq:gar2}
\end{equation}

In the previous argument the value of $0<r<1$ (and thus $s$) was unfixed.
Now $r$ can be optimized to $r=\frac12\left(\sqrt{n^2+6n+1}-n-1\right)$, but, for the sake of simplicity, we take $r=1-\frac2n$.
With this latter choice,
\[\frac{1-r}s\frac1{\vartheta_n}= {\frac {{2}^{1-n}}{n} \left( {\frac {n-2}{n-1}} \right) ^{ n}}\frac1{\vartheta_n} >
  {\frac {{2}^{1-n}}{n} \left( {\frac {1-\frac2n}{\mathrm e}} \right)  }\frac1{\vartheta_n} .
\]

Putting this into  \eqref{eq:gar2}, we obtain  \eqref{eq:gar3}.
\end{proof}

\end{theorem}
\snewpage
From Part I \cite{L1} we may recall:
For $\lambda\in[0,1]$, the convergence radius of $\Theta^{(\lambda)}(x)$ around $x=0$ is
\[{\mathrm C}_\infty^{(\lambda)}=\begin{cases}
2&\text{if }\lambda=\frac12,
\\
\dfrac{2\artanh (1-2\lambda)}{1-2\lambda}=\dfrac{\log\dfrac{1-\lambda}{\lambda}}{1-2\lambda}&\text{if }\lambda\in(0,1)
\setminus\{\frac12\},\\
+\infty&\text{if }\lambda\in\{0,1\}.
\end{cases}\]
This is a strictly convex, nonnegative function in $\lambda\in(0,1)$,
symmetric for $\lambda\mapsto1-\lambda$;  its minimum is ${\mathrm C}_\infty^{(1/2)}=2$.
In particular, in $\lambda\in[0,1]$, it yields a $[2,+\infty]$-valued strictly convex continuous function.

For $d\geq3$, let $\Lambda_d$ be the number such that $0<\Lambda_d\leq\frac12$ and
${\mathrm C}_\infty^{(\Lambda_d )}={\mathrm C}^{\{\{d\}\}}_{\mathbb R}$.
We know that $\Lambda_d$ is increasing in $d$ but $\leq\frac12$.
Then
\[\Lambda_d\geq \Lambda_3=0.0588740902\ldots>\frac1{17}.\]
For a better appreciation, asymptotically,
\begin{lemma} 
\plabel{lem:appr}
For $d\geq25$,
\begin{equation}
\Lambda_d\geq \frac12\left(1-\sqrt{\frac{3}{\lfloor\sqrt{d}\rfloor-1}}\right).
\plabel{eq:emo}
\end{equation}
\begin{proof}
Let $\tilde\Lambda_d$ denote the RHS of \eqref{eq:emo}.
It is sufficient to prove that ${\mathrm C}_\infty^{(\tilde\Lambda_d )}\geq{\mathrm C}^{\{\{d\}\}}_{\mathbb R}$.
In turn, by \eqref{eq:ram},
 for that, it is sufficient to prove that ${\mathrm C}_\infty^{(\tilde\Lambda_d )}\geq2+\frac2{\lfloor\sqrt d\rfloor-1}$.
On the other hand, it is not hard to prove that
${\mathrm C}_\infty^{(\lambda)}\geq2+\frac83(\lambda-1/2)^2$, leading to the desired result.
\end{proof}
\end{lemma}
In general, however, the estimate $\Lambda_d>\frac1{17}$  may be convenient.
Now we can state
\begin{theorem}
\plabel{thm:fin}
Assume that $\mathfrak A$ is Banach algebra of finite real dimension $n$,
 and $\phi$ is an $\mathfrak A$ valued ordered measure such that
\begin{equation}
\int|\phi|<\frac{2}{1-\left( {\dfrac {{2}^{-2-n}}{n}
\left( {\dfrac {1-\frac2n}{\mathrm e}} \right)  }\dfrac{\Lambda_n}{\vartheta_n}\right)}
.
\plabel{eq:dana}
\end{equation}
Then $\psi$ is $M$-controlled, in particular, its Magnus expansion convergent, and the logarithmic Magnus formula holds.
\begin{proof}
Assume that $\phi$ is subject to \eqref{eq:dana}.
Assume $\int|\phi|\geq2$.
Let us decompose $\phi=\phi_1\boldsymbol.\phi_2\boldsymbol. \phi_3$, where
$\int|\phi_1|=1$ and $\int|\phi_2|=2$.
Let $\lambda\in[\Lambda_n,1-\Lambda_n]$.
By Theorem \ref{thm:gain2},
\[\Theta^{(\lambda)}(\phi_1)\leq \Theta^{(\lambda)}( 1)
-\underbrace{\frac{1}2\left( {\frac {{2}^{1-n}}{n}
\left( {\frac {1-\frac2n}{\mathrm e}} \right)  }\frac{\Lambda_n}{\vartheta_n}\right)}_{r:=}
.\]
That means that $\phi_1$ has delay
\[1- \left(\Theta^{(\lambda)} \right)^{-1}(\Theta^{(\lambda)}( 1)-r)\geq
1- \left(\Theta^{(1/2)} \right)^{-1}(\Theta^{(1/2)}( 1)-r)\equiv\frac{r}{4-r}.
\]
(See the ``Delay estimate reduction principle'' of Part I.)
The same applies to $\phi_2$, consequently $\phi_1\cdot\phi_2$ collects delay at least $\frac{2r}{4-r}$.
Thus the resolvent expansion will be convergent (thus the resolvent $\mathcal R^{(\lambda)}(\Rexp \phi)$ exists) as long as
$ \int|\phi_3|<\frac{2r}{4-r}$, i. e. if $\int|\phi |<2+\frac{2r}{4-r}=\frac{2}{1-\frac14r}$.
But that was exactly our assumption.
$\mathcal R^{(\lambda)}(\Rexp \phi)$ also exists if $\int|\phi|<2$.
More generally, the statement about the existence of the resolvent for
$\lambda\in[\Lambda_n,1-\Lambda_n]$
also holds if $\phi$ is replaced by $t\cdot \phi$ where $t$ is from the closed complex unit disk.

Now, we claim, in general,
\begin{equation}
\frac{2}{1-\frac14r}\leq \mathrm C_{\mathbb R}^{\{\{n\}\}}.
\plabel{eq:anto}
\end{equation}
Indeed, otherwise, $\frac{2}{1-\frac14r}> \mathrm C_{\mathbb R}^{\{\{n\}\}}$, and we could chose an $\delta>0$ such that
\[\mathrm C_{\mathbb R}^{\{\{n\}\}}<\frac{2}{1-\frac14r(1-\delta)}.\]
We also know that
 \[\mathrm C_{\mathbb R}^{\{\{n\}\}}=\mathrm C_\infty^{(\Lambda_n)} <\mathrm C_\infty^{(\Lambda_n(1-\delta))}.\]
Finally, we could choose a counterexample $\psi$ to the convergence of the Magnus expansion such that
\begin{equation}\int|\psi|<\frac{2}{1-\frac14r(1-\delta)}\plabel{eq:hora}\end{equation}
and
\begin{equation}\int|\psi|<\mathrm C_\infty^{(\Lambda_n(1-\delta))} .\plabel{eq:xora}\end{equation}
Repeating the arguments of the previous paragraph,
\eqref{eq:hora} causes  $\mathcal R^{(\lambda)}(\Rexp (t\cdot\psi))$ to exist for
$\lambda\in[\Lambda_n(1-\delta),1-\Lambda_n(1-\delta)]$ and $t\in\Dbar(0,1)$;
while \eqref{eq:xora} causes $\mathcal R^{(\lambda)}(\Rexp (t\cdot\psi))$ to exist for $\lambda\in[0,\Lambda_n(1-\delta))$
and  $\lambda\in(1-\Lambda_n(1-\delta),1]$.
Ultimately, $\psi$ is $M$-controlled, and its Magnus expansion is convergent; this is  a contradiction.

Returning to the main statement now, \eqref{eq:anto} implies
\[\int|\phi|\leq \mathrm C_{\mathbb R}^{\{\{n\}\}}\equiv \mathrm C_\infty^{(\Lambda_n)}.\]
This, in turn, implies that  $\mathcal R^{(\lambda)}(\Rexp (t\cdot\phi))$ to exist for $\lambda\in[0,\Lambda_n )$
and  $\lambda\in(1-\Lambda_n ,1]$.
Ultimately, $\phi$ is $M$-controlled.
\end{proof}
\end{theorem}
If we ignore $M$-controlledness and the logarithmic issues, and we concentrate on the
convergence of the Magnus expansion, then the previous theorem says
\begin{equation}
\mathrm C_{\mathbb R}^{\{\{n\}\}}\geq
\frac{2}{1-\left( {\dfrac {{2}^{-2-n}}{n}
\left( {\dfrac {1-\frac2n}{\mathrm e}} \right)  }\dfrac{\Lambda_n}{\vartheta_n}\right)}.
\plabel{eq:subdana}
\end{equation}
Note, however, that there is a huge gap between \eqref{eq:ram} and \eqref{eq:subdana}, which would be interesting to close.

One can improve upon the cumulative convergence radius
$\mathrm C_2=2.89847930\ldots$ of the BCH expansion by similar arguments.

\snewpage

\end{document}